\documentclass[a4j,11pt]{article}

\usepackage{ulem}

\usepackage{amsthm}
\usepackage{amsmath}
\usepackage{amssymb}
\usepackage{amsfonts}
\usepackage{ascmac}
\usepackage{wrapfig}
\usepackage{graphicx,color}
\usepackage{lscape}
\usepackage{comment}
\usepackage{mathrsfs}
\usepackage{fancybox} 

\theoremstyle{plain}
\newtheorem{RC}{R}
\newtheorem{thm}{Theorem}
\newtheorem{lem}{Lemma}
\newtheorem{remark}{Remark}[section]
\newtheorem{coro}{Corollary}

\begin{document}
\title{Asymptotic Moments Matching to Uniformly Minimum Variance Unbiased Estimation under Ewens Sampling Formula}

\author{Masayo Y. Hirose\thanks{Institute of Mathematics for Industry, Kyushu University, masayo@imi.kyushu-u.ac.jp} \and Shuhei Mano\thanks{Department of Statistical Inference and Mathematics, The Institute of Statistical Mathematics, smano@ism.ac.jp}}

\date{}

\maketitle
\begin{abstract}
The Ewens sampling formula is a distribution related to the random partition of a positive integer. In this study, we investigate the issue of non-existence solutions in parameter estimation under the distribution. As a result, the first and second moments matching estimators to the uniformly minimum variance unbiased estimator are derived using the Ewens sampling formula in asymptotic sense. A Monte Carlo simulation study is performed to evaluate the efficiency of the resulting estimators.

\end{abstract}

{\bf Keywords:} Ewens Distribution, Higher-Order Efficiency, Population Unique.

\section{Introduction}

Ewens (1972) provided the Ewens sampling formula, which is a law of the partition of positive integers into components comprising non-negative integers in the context of genetics. Antoniak (1974) derived the Ewens sampling formula in the context of Bayesian statistics as a partition induced by a sample from the Dirichlet process. It has been applied to several research fields in ecology, disclosure control and so on. Let $S\in \mathbb S_N$ satisfy

$$\mathbb S_N=\{S\equiv (S_1,\ldots, S_N): S_i\geq 0, \sum_{i=1}^N i S_i=N, i\in \{1,\ldots,N\}\}.$$ 

Elements $S_i$ of the $N$-dimensional vector $S$ are random variables and denote the number of types of $i$ times appear. This is known as the frequency of the frequencies (Good, 1953). Next, we define $K_N=\sum_{i=1}^N S_i$.  The number of types to appear is denoted, i.e., the length of the random partition. For instance, in ecology, $S_i$ and $K_N$ denote the number of species that appear $i$ times and the number of different species to occur, respectively. For a positive integer $N$, it is expressed using the parameter $\theta>0$ as follows:
\begin{align}
P(S=(s_1,\ldots,s_N)\in \mathbb S_N)=\frac{N!\theta^{k}}{\theta^{[N]}}\prod_{j=1}^N\frac{1}{j^{s_j}s_j!},\label{e1.210413}
\end{align}
where $k=\sum_{i=1}^N s_i$ and $\theta^{[N]}=\prod_{j=1}^N(\theta+j-1)$. 

\noindent The parameter $\theta$ controls the ``diversity'': $k \to N$ as $ \theta \to \infty$, and $k \to 1$ as $\theta \to 0$. It is of particular interest to assess $S_1$. For instance, $S_1$ denotes the number of singleton species in the population and a population unique. The latter is explained in Section 4 as an example. For more details regarding this formula, readers may refer to Tavare and Ewens (1997), Crane (2016), and Mano (2018).

The distribution of $K_N=\sum_{i=1}^N S_i$, given $N$, can be obtained as follows:
\begin{align}
P(K_N=k)=&\frac{\theta^{k}}{\theta^{[N]}}s(N,k), \label{Un.dist}
\end{align}
where $s(N,k)$ is the unsigned Stirling numbers of the first kind satisfying $\theta^{[N]}=\sum_{k=1}^N s(N,k) \theta^{k}$ for non-negative integers $k$ and $N$ such that $1\leq k\leq N<\infty$.

We let $R_i=E[S_i]$ be the expectation of $S_i$ with respect to this model. Its explicit formula is obtained as a function of $\theta$ in Watterson (1974) as follows:
\begin{align*}
R_{i}=R_{i}(\theta)=\frac{\theta}{i}\prod_{j=1}^i \frac{N-j+1}{\theta+N-j}.
\end{align*}

Hereinafter, we focus primarily on the inference for $R_i$ in the population from the sample data. When the Ewens sampling formula is supposed to be a population model, $N$ denotes the population size, whereas $N$ is replaced with the sample size $n$ in formula \eqref{e1.210413} when the Ewens sampling formula indicates the sampling distribution. This is reasonable because the Ewens sampling formula has a partition structure (Kingman, 1978) which the distribution of exchangeable random partitions coincides with that of any subsampling partitions with sample size $n$ from population size $N$ for all $n<N$. Therefore, in practice, for sample size $n$, an estimator of $R_i$ can be obtained by replacing the unknown parameter $\theta$ in the formula of $R_i$ with its consistent estimator.

In particular, the maximum likelihood estimator $\hat \theta_{ML}$ is widely used and is obtained as follows:
\begin{align*}
\hat \theta_{ML}=\arg \max_{\theta \in \mathbb R_{>0}}l(\theta),
\end{align*}
where $\mathbb R_{>0}$ is the parameter space of $\theta$, and the logarithm of the likelihood $l(\theta)$ is such that
\begin{align*}
l(\theta)=K_n \log(\theta)-\sum_{j=1}^n\log (\theta+j-1)+const. 
\end{align*}

Tavare and Ewens (1997) reported, ``This estimator is biased, but the bias decreases as $n$ increases.'' In addition, if $R_i$ is U-estimable, it may be difficult to obtain a uniformly minimum variance unbiased estimator (UMVUE).  To the best of our knowledge, asymptotic moments matching estimators to a uniformly minimum variance unbiased estimator of $R_i$ have not been investigated, whereas $E[K_n]$ is an exact UMVUE. We denote this estimator the asymptotic UMVUE hereinafter. Moreover, even if we can identify such an estimator, it may yield unrealistic negative estimates of non-negative $R_i$.

First, we address the construction of the asymptotic UMVUE of $R_i$ with a small positive integer $i$, up to the second moment matching in second-order asymptotic sense. In this study, the second-order indicates the order $O((\log n)^{-1})$ for large $n$. We also avoid severe issues such that the maximum likelihood solution could not exist (see Section 2). In addition, the precision of the estimation can be improved further. Second, we construct a higher-order asymptotic UMVUE by matching the first and second moments.

For the second purpose, we use two types of bias correction methods: additive bias correction and the adjusted maximum likelihood method. The adjusted maximum likelihood is multiplied a nonrandom adjustment factor by the likelihood. It was developed by Lahiri and Li (2010), Li and Lahiri (2010), and Yoshimori and Lahiri (2014) to avoid zero estimates of dispersion parameters in a linear mixed model, particularly in the research fields of small-area estimation. Hirose and Lahiri (2018) achieved a second-order asymptotic unbiasedness of several important re-parameterized estimators by suggesting a new adjusted maximum likelihood method. Hirose and Mano (2021) constructed a general framework using differential geometry to achieve second-order unbiasedness and applied the methodology to a general model with multi-dimensional parameters.

One may consider these previous results to be applicable because the Ewens distribution is an exponential distribution family. Nevertheless, these results are insufficient to obtain the asymptotic UMVUE up to the fourth-order asymptotic sense, unlike the studies of Hirose and Lahiri (2018) and Hirose and Mano (2021). It is note that, in this study, the second (or fourth) order denotes the order of $O((\log n)^{-1})$ (or $O((\log n)^{-2})$) for large $n$, whereas the second-order denotes the $O(n^{-1})$ in the study of Hirose and Mano (2021), considering the differences in the Fisher information order.

As a result, it is sufficient to obtain one common estimate of the parameter for achieving the second purpose even when $R_i$ and $R_j$ are to be estimated simultaneously for $i\neq j$. Furthermore, we demonstrate the higher-order asymptotic results based on easier proofs, owing to the functional form $R_i$ of $\theta$, a property of the exponential distribution family and a relationship between the parameter $\theta$ and natural parameter $\xi$. For more details, see Section 3.2 and Appendix C.3.

The remainder of this paper is organized as follows: In Section 2, we introduce the existing estimator of $R_i$ and modify it to be the asymptotic UMVUE, up to the second-order, to achieve the first purpose. The problem of non-existence of estimates is avoided in this section. To address the second problem, higher-order asymptotic unbiasedness is discussed in Section 3. In this section, we suggest two estimators using two bias-correction methods. This methodology can be applied to practical issues. Subsequently, in Section 4, we present an example where our methodology is applied to estimate the number of population uniques to disclosure control. The simulation study is described in Section 5. Herein, it is assumed that $\theta$ and $i$ are bounded for large $n$. In addition, we assume that the sampling design is simple random sampling without replacement. All technical proofs are provided in the appendix.

\section{Maximum Likelihood Estimation of Parameter}

As mentioned earlier, a typical method to estimate $R_i$ is to replace the $\theta$ of $R_i$ with its maximum likelihood estimator. Its first and second derivatives are expressed as follows:
\begin{align}
\frac{\partial l(\theta)}{\partial \theta}\equiv \partial_\theta l(\theta)=&\frac{K_n}{\theta}-\sum_{j=1}^n \frac{1}{\theta+j-1},\label{e1.1103}\\
\frac{\partial^2 l(\theta)}{\partial \theta^2}\equiv \partial_\theta^2 l(\theta)=&-\frac{K_n}{\theta^2}+\sum_{j=1}^n \frac{1}{(\theta+j-1)^2}.\label{e2.210410}
\end{align}

\noindent From \eqref{e1.1103}, the maximum likelihood estimator $\hat \theta_{ML}$ can be expressed as the root of 
$$\frac{K_n}{\theta}-\sum_{j=1}^n \frac{1}{\theta+j-1}=0. $$ 
We let $\hat R_{i}^{(N)}=R_{i}(\hat \theta_{ML})$ and refer to it as a {\it naive estimator}.

The continuous mapping theorem may provide the consistency of the naive estimator; to the best of our knowledge, the properties of asymptotic moments have not been investigated hitherto.

The following practical issue occurs when attempting to obtain an asymptotic UMVUE: cases $K_n \in \{1, n\}$ provide each likelihood as a strict monotone function of $\theta$. Additionally, it is shown from the first-order derivative of the likelihood function of the natural parameter 
$\xi=\log \theta \in \Xi \subset \mathbb R$ on the Ewens distribution.

\begin{align*}
\partial_\xi l(\xi)=&K_n-n+\sum_{j=2}^n \frac{j-1}{e^{\xi}+j-1}.
\end{align*}
\noindent It is note that the equation above is rewritten from \eqref{e1.1103} with the natural parameter $\xi$.

The following lemma is established to assess such probability.
\begin{lem}
 \label{lem2}
Under the regularity condition R1, we have the following for large $n$:
\begin{align*}
P(\mathscr K)=o((\log n)^{-2}), 
\end{align*}
where the set $\mathscr K=\{K_n\in\{1,n\}\}$. 
\end{lem}
\noindent The regularity condition and proof are provided in Appendices A and C.1, respectively.

In set $\mathscr K$, the lemma implies the problem of non-existing maximum likelihood solutions in 
$\Theta\subset \mathbb R_{>0}$ with extremely low but non-zero probability. 
In addition, $\hat R_{i}^{(N)}$ may be absent.

Therefore, the estimator must be modified to obtain an asymptotic UMVUE by considering such cases. Hence, we define set 
$\mathscr{S}=\{\tilde \theta: \tilde \theta \in (0,C_+]\}$ with a large positive finite value $C_+$, which does not depend on $n$. Subsequently, we let
\begin{eqnarray*}
\hat R_{i}^{(NM)}\equiv\left\{
\begin{array}{ll}
0 &(K_n=1)\\
R_{i}(\hat \theta_{ML}) &(\hat \theta_{ML} \in \mathscr{S})\\
R_{i}(C_+)& (\hat \theta_{ML} \notin \mathscr{S}\cup \{0\}). 
\end{array}
\right.
\end{eqnarray*} 
\noindent It is note that $R_i (C_+)$ can be adopted in case $\{K_n=n\}$ because $\{K_n=n\}\subset  \mathscr{S}^c\cap \{0\}^c$. In addition, $\hat \theta_{ML}>0$ when $1<K_n<n$ because $\xi$ exists such that $\partial_\xi l(\xi)=0$.

Next, another lemma is established, the proof of which is provided in Appendix C.2. 
\begin{lem}
 \label{lem2.5}
Under regularity conditions R1 and R2, we have the following for large $n$:
\begin{align*}
P(\mathscr S^c)=o((\log n)^{-2}).
\end{align*}
\end{lem}

\noindent Therefore, the estimator $\hat R_{i}^{(NM)}$ is the function of the complete sufficient statistic $K_n$ of $\theta$ in cases $\hat \theta_{ML}\in \mathscr S$. From Lemmas \ref{lem2} and \ref{lem2.5}, Theorem \ref{thm0} shows that estimator $\hat R_i^{(NM)}$ is the second-order asymptotic UMVUE for large $n$.

\begin{thm}
\label{thm0}
Under regularity condition R1, the following holds:
\begin{align*}
(i) &E[\hat R_{i}^{(NM)}-R_i(\theta)]=O((\log n)^{-2}),\\
(ii) &E[\{\hat R_{i}^{(NM)}-R_i(\theta)\}^2]=\frac{R_{i}^2}{\sum_{j=2}^n\frac{\theta (j-1)}{(\theta+j-1)^2}}+o((\log n)^{-2}).
\end{align*}
\end{thm} 
\noindent This is shown in Appendix B.1.

Next, we provide a remark.

\begin{remark}
\label{rem1}
The length of the random partition $K_N$ in the population is also of interest to infer. It holds that
\begin{align*}
\eta(\theta)\equiv E[K_N]&=\sum_{j=1}^N\frac{\theta}{\theta+j-1}. 
\end{align*} 
We do not address improving this estimator because the estimator $\eta(\hat \theta_{ML})$ becomes the exact UMVUE of $\eta(\theta)$. The result is a well-known result and can be an example of Corollary 3.13 in Hirose and Mano (2021).

\end{remark}

\section{Higher-Order Asymptotic UMVUE}

\subsection{General bias-corrected estimator}

Theorem \ref{thm0} shows that the bias is of the order $O((\log n)^{-2})$ for large $n$. In this section, we address the construction of two types of asymptotic UMVUEs for matching the first and second moments in the fourth-order asymptotic sense for large $n$. In Sections 3.2 and 3.3, we present the results using additive bias correction and the adjusted maximum likelihood method, respectively.

Let $\hat R_{i}^{(BC)}$ denote the general bias-corrected estimator of $R_i$ while considering set $\mathscr K$.
\begin{eqnarray}
\hat R_{i}^{(BC)}\equiv\left\{
\begin{array}{ll}
0 &(K_n=1)\\
\hat R_{i}(\hat \theta) &(\hat \theta \in \mathscr{S})\\
R_{i}(C_+)& (\hat \theta \notin \mathscr{S}\cup \{0\}),
\end{array}
\right.\label{e1.210127}
\end{eqnarray} 
provided $\hat \theta>0$, where $\hat R_{i}(\hat \theta)$ and $\hat \theta$ denote consistent estimators of $R_i$ and $\theta$ for large $n$, respectively.
We set $\hat \theta$ as a small positive value when the solution of $\theta$ is negative, except for $K_n=1$. However, it holds that $\hat \theta_{ML}>0$ when $1<K_n<n$.

\subsection{Additive bias correction for the higher-order asymptotic inference}

One may consider using the additive bias correction method to reduce bias.  Let $\hat R_{i}^{(BC1)}$ be the term $\hat R_{i}(\hat \theta)$ in \eqref{e1.210127} is replaced with $\hat R_{i}^{(N)}-B_i(\hat \theta_{ML})$, where 
$$B_i(\theta)=\frac{\sum_{j=2}^n\frac{(j-1)}{(\theta+j-1)^3}}{\left\{\sum_{j=2}^n\frac{(j-1)}{(\theta+j-1)^2}\right\}^2}
R_i.$$

Then, the theorem establishes that $\hat R_{i}^{(BC1)}$ achieves the fourth-order unbiasedness for large $n$, while maintaining the asymptotic efficiency.

\begin{thm}
\label{thm2}
Under regularity condition R1 for large $n$, the following holds:
\begin{align*}
(i) &E[\hat R_{i}^{(BC1)}-R_{i}(\theta)]=o((\log n)^{-2}),\\
(ii) &E[\{\hat R_{i}^{(BC1)}-R_{i}(\theta)\}^2]=\frac{R_{i}^2}{\sum_{j=2}^n\frac{\theta (j-1)}{(\theta+j-1)^2}}+o((\log n)^{-2}). 
\end{align*}
\end{thm} 
\noindent The proof is provided in Appendix B.1.

In addition, the estimator $\hat R_{i}^{(BC1)}$ is a function of the complete sufficient statistic $K_n$ of $\theta$ in cases $\hat \theta_{ML}\in \mathscr S$. Therefore Lemmas \ref{lem2}, \ref{lem2.5}, and Theorem \ref{thm2} prove that it is the asymptotic UMVUE, up to the fourth-order.

\subsection{Another possible bias correction: adjusted maximum likelihood method}

Alternatively, bias can also be reduced using the adjusted maximum likelihood method. This method has been developed by Hirose and Lahiri (2018) and Hirose and Mano (2021) for bias correction after re-parameterization. Firth (1993) suggested a similar bias reduction method for $\hat \theta_{ML}$ via second-order asymptotic expansion using a score function. In this section, unlike their methods, we present the derivation of the fourth-order asymptotic UMVUE using the higher-order asymptotic expansion.

We define the general adjusted maximum likelihood estimator of $\theta$ as follows:
\begin{align}
\hat \theta_{GA}=\arg \max_{\theta>0} l_{ad}(\theta),\label{e2.210127}
\end{align}
where $l_{ad}(\theta)=l(\theta)+\tilde l_{ad}(\theta).$

\noindent In addition, we denote $e^{\tilde l_{ad}}$ as the adjustment factor. For example, the maximum likelihood estimator $\hat \theta_{ML}$ is obtained when $\tilde l_{ad}(\theta) \propto C$ is adopted, where $C$ is a constant value that does not depend on $\theta$.

Next, we let $\hat R_{i}^{(BCA)}$ be an estimator of $R_i$, where $\hat R_{i}$ in \eqref{e1.210127} is replaced with $R_{i}(\hat \theta_{GA})$. 
Theorem \ref{thm3} is presented to show its property of asymptotic moments, of which the proof is shown in Appendix B.2.

\begin{thm}
\label{thm3}
For large $n$ under regularity conditions R1 and R2,
\begin{align*}
(i) E[\hat R_{i}^{(BCA)}-R_{i}(\theta)]=\frac{R_i}{g_\xi}&
\left[\partial_\xi^{(1)} \tilde l_{ad}\left(1+\frac{\partial_\xi^{(2)}\tilde l_{ad}}{g_{\xi}}\right)-\frac{\partial_\xi^{(2)} \tilde l_{ad}\partial_\xi g_\xi}{2g_{\xi}^2}\right]\notag\\
&+\frac{R_i}{2g_\xi^2}\left(g_\xi-\partial_\xi g_\xi-\partial_\xi^{(2)} \tilde l_{ad}+\partial_\xi^{(3)} \tilde l_{ad}\right)+o((\log n)^{-2});\notag\\
(ii) E[\{\hat R_{i}^{(BCA)}-R_{i}(\theta)\}^2]=&\frac{R_{i}^2}{g_\xi}\left(1+\frac{(\partial_\xi \tilde l_{ad}(\xi))^2
+2\partial_\xi^{(2)} \tilde l_{ad}(\xi)}{g_\xi}\right)+o((\log n)^{-2}),
\end{align*}
where $\xi=\log \theta$ and $g_\xi=\sum_{j=2}^n\frac{\theta (j-1)}{(\theta+j-1)^2}$.
\end{thm}

The theorem above implies that the following condition of the adjustment factor is required to eliminate the fourth-order asymptotic bias without sacrificing the asymptotic efficiency.
\begin{align}
\partial_\xi \tilde l_{ad}(\xi)\left(1+\frac{\partial_\xi^{(2)}\tilde l_{ad}}{g_{\xi}}\right)-\frac{\partial_\xi^{(2)} \tilde l_{ad}\partial_\xi g_\xi}{2g_{\xi}^2}&+\frac{1}{2g_\xi}\left(g_\xi-\partial_\xi g_\xi-\partial_\xi^{(2)}\tilde l_{ad}(\xi)+\partial_\xi^{(3)}\tilde l_{ad}(\xi)\right)\notag\\
&=o((\log n)^{-1}),\label{e1.210129}\\
(\partial_\xi \tilde l_{ad}(\xi))^2+2\partial_\xi^{(2)} \tilde l_{ad}(\xi)&=o(1). \label{e3.210124}
\end{align}

\noindent To obtain an adjustment factor satisfying \eqref{e1.210129} and \eqref{e3.210124}, we therefore restrict the class of adjustment factor to the following for large $n$ with $j=1,2,3$:
$$\partial_\theta^{(j)} \tilde l_{ad}(\theta)=o(1).$$

\noindent We then find that the resulting specific adjustment factor $e^{\tilde l_{ad}}$ satisfies
\begin{align}
\partial_\theta  \tilde l_{ad}(\theta)=&-\frac{1}{\theta}\frac{g-\partial_\xi g}{2g}=-\frac{\sum_{j=2}^n\frac{(j-1)}{(\theta+j-1)^3}}{\sum_{j=2}^n\frac{(j-1)}{(\theta+j-1)^2}}.\label{e3.210410}
\end{align}
Additionally, it holds that $\partial_\xi^{(j)}  \tilde l_{ad}(\xi)=O((\log n)^{-1})$ with $j=1,2,3$ from \eqref{e1.1.18} given in Appendix A.

Subsequently, we let $\hat \theta_A$ be the above-mentioned adjusted maximum likelihood estimator. It can be obtained as the root of the following equation:
$$\frac{K_n}{\theta}-\sum_{j=1}^n\frac{1}{\theta+j-1}-\frac{\sum_{j=2}^n\frac{(j-1)}{(\theta+j-1)^3}}{\sum_{j=2}^n\frac{(j-1)}{(\theta+j-1)^2}}=0.$$ 
In addition, we define the estimator $\hat R_i^{(BC2)}$, which substitutes $\hat R_i$ in \eqref{e1.210127} with $R_i(\hat \theta_{A})$.
Corollary \ref{coro1} summarizes the fourth-order asymptotic properties of its first and second moments.

\begin{coro}
\label{coro1}
Under regularity condition R1 for large $n$, the following holds:
\begin{align*}
(i) &E[\hat R_{i}^{(BC2)}-R_{i}(\theta)]=o((\log n)^{-2});\\
(ii) &E[\{\hat R_{i}^{(BC2)}-R_{i}(\theta)\}^2]=\frac{R_{i}^2}{\sum_{j=2}^n\frac{\theta (j-1)}{(\theta+j-1)^2}}+o((\log n)^{-2}).
\end{align*}
\end{coro}

Corollary \ref{coro1} shows that our specific adjustment factor contributes to the disappearance of the fourth-order asymptotic bias without sacrificing asymptotic efficiency. As it is for $\hat R_{i}^{(BC1)}$, the estimator is a function of the complete sufficient statistic $K_n$ of $\theta$ in cases $\hat \theta_{A}\in \mathscr S$. Therefore, it is also an asymptotic UMVUE up to the fourth-order, based on Lemmas \ref{lem2} and \ref{lem2.5}.

Next, we provide some remarks.

\begin{remark}
The estimator $\hat R_{i}^{(BC1)}$ makes a downward correction by the term $\hat B_i$ because the bias-corrected term $\hat B_i$ is positive almost surely. Nevertheless, the term $B_i$ has a slightly complex functional structure, which may result in unrealistic negative estimates of $\hat R_{i}^{(BC1)}$. By contrast, $\hat R_{i}^{(BC2)}$ maintains a simple function structure and ensures that it is in the range of $R_i(\theta)$ simultaneously. 
\end{remark}

\begin{remark} 
It is noteworthy that the logarithm of the adjustment factor $\tilde l_{ad}(\theta)$ does not depend on $i$ for obtaining $\hat \theta_A$. In other words, $\hat \theta_A$ can be used as one common estimate of $\theta$ even in the simultaneous higher-order asymptotic inferences of $\hat R_{i}$ and $\hat R_{j}$ ($i\neq j$), while maintaining the function form of $R_i$. This may significantly reduce the computer burden in simultaneous inferences.
\end{remark}

\begin{remark} 
The bias correction for $R_i(\theta)$ using our adjusted maximum likelihood method corresponds to the bias correction for $\hat \theta_{ML}$, owing to the results \eqref{e1.1.18} associated with the orders of $\partial_\theta R_{i}$ and $\partial_\theta^{2} R_{i}$. 
One may recall that Firth (1993) also suggested a bias reduction method for $\hat \theta_{ML}$. However, we expanded the higher-order asymptotic expansion and derived a fourth-order asymptotic UMVUE. 
\end{remark}

\begin{remark} 
The Ewens model enables the asymptotic theoretical result of $R_i$ to be constructed easily, owing to the properties of the exponential family distribution and the relationship $\theta=e^{\xi}$, where $\xi$ is the natural parameter. For more details, see Appendix C.2. In general, a more complex proof may be required for the fourth-order asymptotic expansion. 
\end{remark}

\begin{remark} 
Our adjustment factor obtained from \eqref{e3.210410} coincides with one example of Corollary 3.9 in Hirose and Mano (2021) for one-flat manifold, although the asymptotic orders of the Fisher information and expansion are different from those of our study.
\end{remark}

\section{Application to disclosure control: assessed risk of population unique}

In official statistics, data providers often create secondary available tables from microdata for disclosure to users while guaranteeing security. In this case, the microdata are categorized by the attribution of individuals in the cells of the table.  To protect personal information, the risk of individual identification must be assessed from such a table. This risk is referred to as the microdata disclosure risk.

In this example, $S_i$ represents the number of cells on the population, in which the number of individuals is $i$, whereas $K_N$ denotes the number of non-empty cells for population size $N$. 
Sibuya (1993) named the $S_i$ as {\it size indices}.

The indices are used to assess the disclosure risk. For instance, when $S_i$ is large with a small $i$, it is interpretable that the table has a high disclosure risk. In practice, $S_i$ is estimated using $R_i$. In particular, $S_1$ is especially of interest to assess and is referred to as the number of {\it population uniques}. By contrast, the cells in which individuals are unique in the sample is known as {\it sample unique}. Additionally, the risk of ``population and sample unique'' is assessed through $f\times R_1$, where $f$ is a known sampling ratio. In such a case, $S_1$ should also be estimated because the number of sample uniques can be observed.

For the inference, super-population models are often used (Bethlehem et al., 1990; Hoshino and Takemura, 1998; Hoshino, 2001). In this study, we assume that the Ewens sampling formula is not only a super-population model, but also a sampling model. As mentioned earlier, this is reasonable because the Ewens sampling formula has a partition structure (Kingman, 1978).

It is clear that the previous methodology can be applied to estimate the number of population uniques as a disclosure risk. 
Theorem \ref{thm0} realizes the asymptotic UMVUE of $R_1$, up to the order of $O((\log n)^{-1})$, as 
\begin{eqnarray*}
\hat R_{1}^{(NM)}\equiv\left\{
\begin{array}{ll}
0 &(K_n=1)\\
R_{1}(\hat \theta_{ML}) &(\hat \theta_{ML} \in \mathscr{S})\\
R_{1}(C_+)& (\hat \theta_{ML} \notin \mathscr{S}\cup \{0\}),
\end{array}
\right.
\end{eqnarray*} 
where $C_+$ denotes a large but finite positive constant.

Next, from Theorems \ref{thm2} and \ref{thm3}, $\hat R_1^{(BC1)}$ and $\hat R_1^{(BC2)}$ become the fourth-order asymptotic UMVUEs of $R_1$. Specifically, $\hat R_i^{(BC1)}$ is expressed as
\begin{eqnarray*}
\hat R_{1}^{(BC1)}\equiv\left\{
\begin{array}{ll}
0 &(K_n=1)\\
R_{1}(\hat \theta_{ML})\left[
1-\frac{\sum_{j=2}^n\frac{(j-1)}{(\hat \theta_{ML}+j-1)^3}
}{\left\{\sum_{j=2}^n\frac{(j-1)}{(\hat \theta_{ML}+j-1)^2}\right\}^2}
\right] &(\hat \theta_{ML} \in \mathscr{S})\\
R_{1}(C_+)& (\hat \theta_{ML} \notin \mathscr{S}\cup \{0\}).
\end{array}
\right.
\end{eqnarray*} 

\noindent An alternative estimator $\hat R_1^{(BC2)}$ provides a simpler formula, as follows:
\begin{eqnarray*}
\hat R_{1}^{(BC2)}\equiv\left\{
\begin{array}{ll}
0 &(K_n=1)\\
R_{1}(\hat \theta_{A}) &(\hat \theta_{A}  \in \mathscr{S})\\
R_{1}(C_+)& (\hat \theta_{A} \notin \mathscr{S}\cup \{0\}),
\end{array}
\right.
\end{eqnarray*} 
provided that $\hat \theta_{A}>0$.

As mentioned in Remark \ref{rem1}, the expectation parameter $\eta$ is used for the inference of the number of non-empty cells $K_N$ to obtain the exact UMVUE.

\section{Monte-Carlo simulation}

We implemented a finite sample simulation study to assess the efficiency of several estimators $\hat R_1$ through Monte-Carlo simulations.

Hence, we considered certain simulation settings such that population size  $N=10^4$, three sample size patterns, i.e., $n\in \{20,10^2,10^3\}$, and $10^4$ replications were generated from the Ewens distribution. Moreover, we set 15 (five values in each of the three patterns P1--P3) patterns of $\theta$ for each sample size $n$ to evaluate the relative effect of the true value of $\theta$ for sample size $n$, as follows: {\bf P1:} $\theta \in \{1,3,5,7,9\}$; {\bf P2:} $\theta \in \{10,30,50,70,90\}$; {\bf P3:} $\theta \in \{100,300,500,700,900\}$. 

\noindent Some cases existed where $n<\theta$, these asymptotic setting of which was not considered to obtain the theoretical result in this study. However, these results were also reported herein.

\noindent Three estimators of $R_1$ were considered for comparison, as follows:
(i) the second-order asymptotic UMVUE $\hat R_1^{(NM)}$, introduced in Section 2; 
(ii) the fourth-order asymptotic UMVUE $\hat R_1^{(BC1)}$, introduced in Section 3.2;
(iii) the fourth-order asymptotic UMVUE $\hat R_1^{(BC2)}$, introduced in Section 3.3.
The estimators of (i)--(iii) are denoted as ``NM,'' ``BC1,'' and ``BC2,'' respectively. In addition, $C_+=10^6$ was adopted.

We first evaluated the estimators $\hat R_i$ using the relative bias and relative root of the mean squared error for the true $R_1$. The relative bias (RB) and relative root of the mean squared error (RRMSE) are defined as
\begin{align*}
RB&\equiv \frac{1}{10^4 \times R_1}\sum_{r=1}^{10^4} (\hat R_1^{(r)}-R_1^{(r)})\times 100,\\ 
RRMSE&\equiv \frac{1}{R_1}\left\{{\frac{1}{10^4}\sum_{r=1}^{10^4} (\hat R_1^{(r)}-R_1^{(r)})^2}\right\}^{1/2}\times 100,
\end{align*}
where an estimate $\hat R_1^{(r)}$ and a true value $R_1^{(r)}$ are constructed using the $r$th replication with $r=1,\ldots, 10^4$.

\begin{figure}[h!]
\begin{center}
 \includegraphics[width=10cm, bb=0 0 1082 1082]{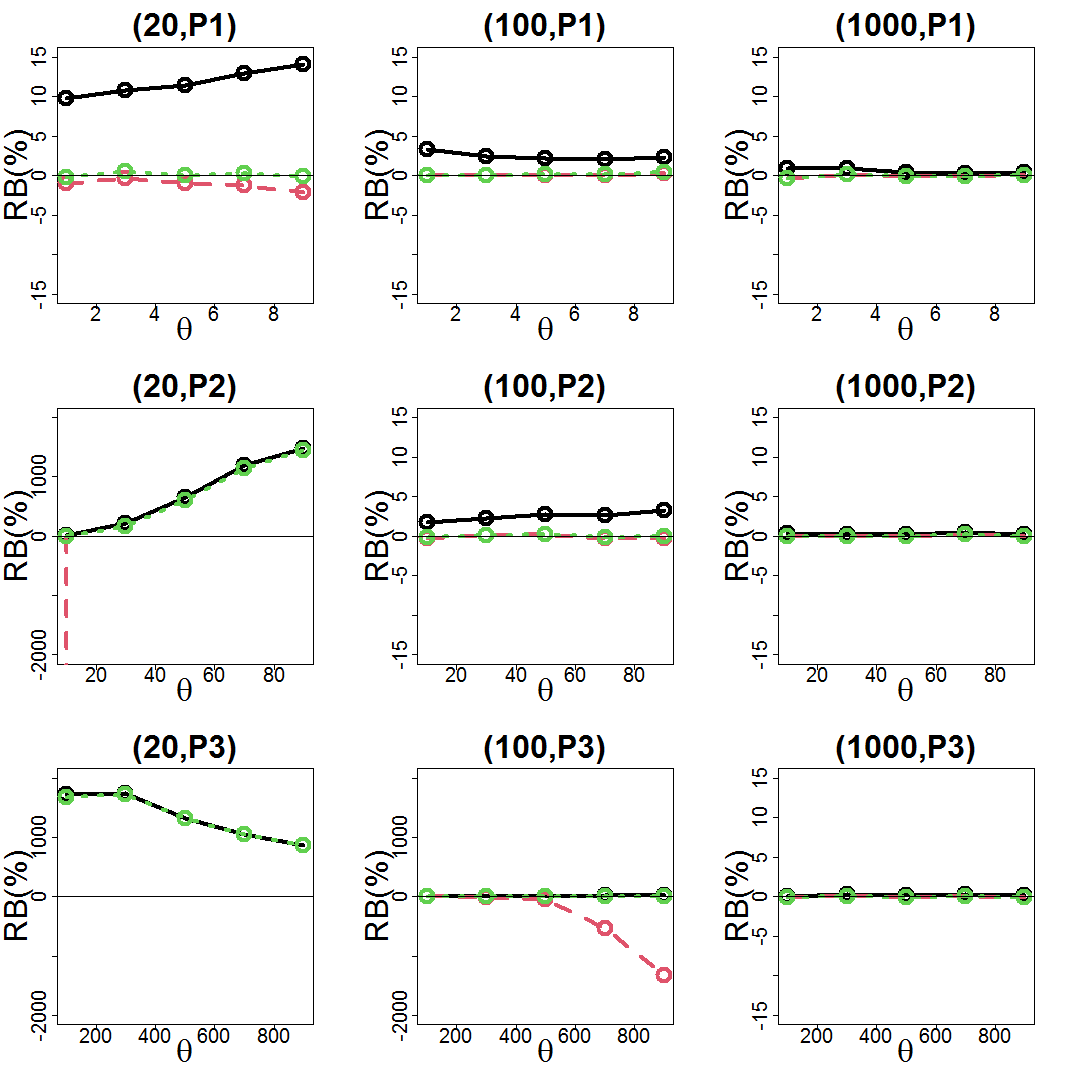}
 \caption{Relative biases (RBs) of three estimators of $R_1$ (NM: solid line; BC1: dashed line; and BC2: 
dotted line) in each combination $(n,P)$, where $P$ is one of three patterns P1--P3 for $\theta$, and sample size $n$ is fixed at $20$ (left), $10^2$ (center), and $10^3$ (right); $x$ axis denotes values of $\theta$.}
\label{RB}
\end{center}
\end{figure}

Figure \ref{RB} shows the relative biases (RB) in nine figures for each combination of $(n,P)$, where $P$ denotes one of three patterns P1--P3 for $\theta$. The right side of the three figures show the results for case $n=10^3$; as shown, all estimators demonstrated similar performance in terms of the relative bias. Meanwhile, in three other figures for cases $(n,P)\in\{(20,P1),(100,P1), (100,P2)\}$, two fourth-order asymptotic UMVUEs that outperformed the second-order asymptotic one are shown. In particular, in $(n,P)=(20,P1)$, $\hat R_1^{(BC2)}$ performed better than $\hat R_1^{(BC1)}$ in terms of the relative bias, as shown in the left of the top figures. Our asymptotic setting in this study did not consider the following simulation settings: $(n,P) \in \{(20,P2),(20,P3),(10^2,P3)\}$. Nonetheless, we also reported these results by changing the scale of the $y$ axis, although some results of $\hat R_{1}^{(BC1)}$ were not appeared because of their considerably low relative biases. Furthermore, these figures show that $\hat R_{1}^{(BC1)}$ can underestimate significantly, whereas the others performed similarly when $n$ is smaller than $\theta$.  This might be caused by the inflation of $\hat B_1$, which may suggest another possibility for the theoretical differences between the fourth-order asymptotic UMVUEs in other asymptotic settings.

\begin{figure}[h!]
\begin{center}
 \includegraphics[width=10cm, bb=0 0 1082 1082]{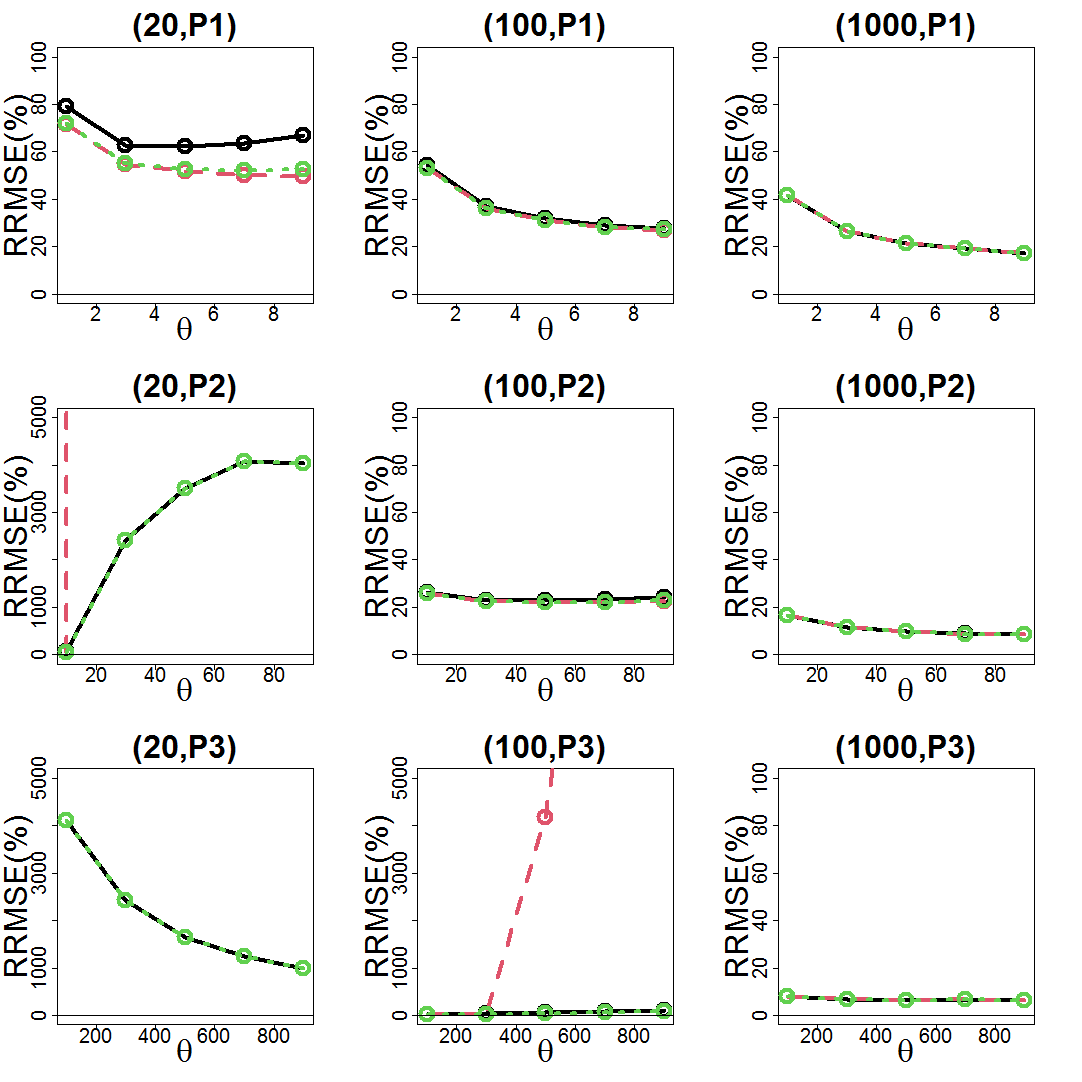}
  \caption{Relative root of mean squared errors (RRMSEs) of three estimators of $R_1$ (NM: solid line; BC1: dashed line; and BC2: dotted line) in each $(n,P)$ combination, 
where $P$ is one of three patterns P1--P3 for $\theta$; sample size $n$ was fixed as $20$ (left), $10^2$ (center), and $10^3$ (right); $x$ axis denotes values of $\theta$.}
\label{RRMSE}
\end{center}
\end{figure}

Next, the relative root of mean squared error (RRMSE) is shown in Figure \ref{RRMSE}, which comprises nine figures for each combination of $(n,P)$.

\noindent The case $(n,P)=(20,P1)$ demonstrates the superiority of the fourth-order asymptotic UMVUEs in terms of the relative root of the mean squared error. In cases $(n,P) \in \{(10^2,P1),(10^2,P2),(10^3,P1),(10^3,P2), (10^3,P3)\}$, we did not observe significant differences among all candidates from the figures. Moreover, as it is for Figure \ref{RB}, we reported three other cases for $n<\theta$ with a scale change for the $y$ axis, and some results of $\hat R_{1}^{(BC1)}$ were not appeared because of their considerably large relative roots of the mean squared errors. Such results might be due to the considerable underestimation of $\hat R_{1}^{(BC1)}$. By contrast, no significant differences were observed between the other two estimators even in such cases.

\begin{figure}[h!]
\begin{center}
 \includegraphics[width=13cm, bb=0 0 2100 700]{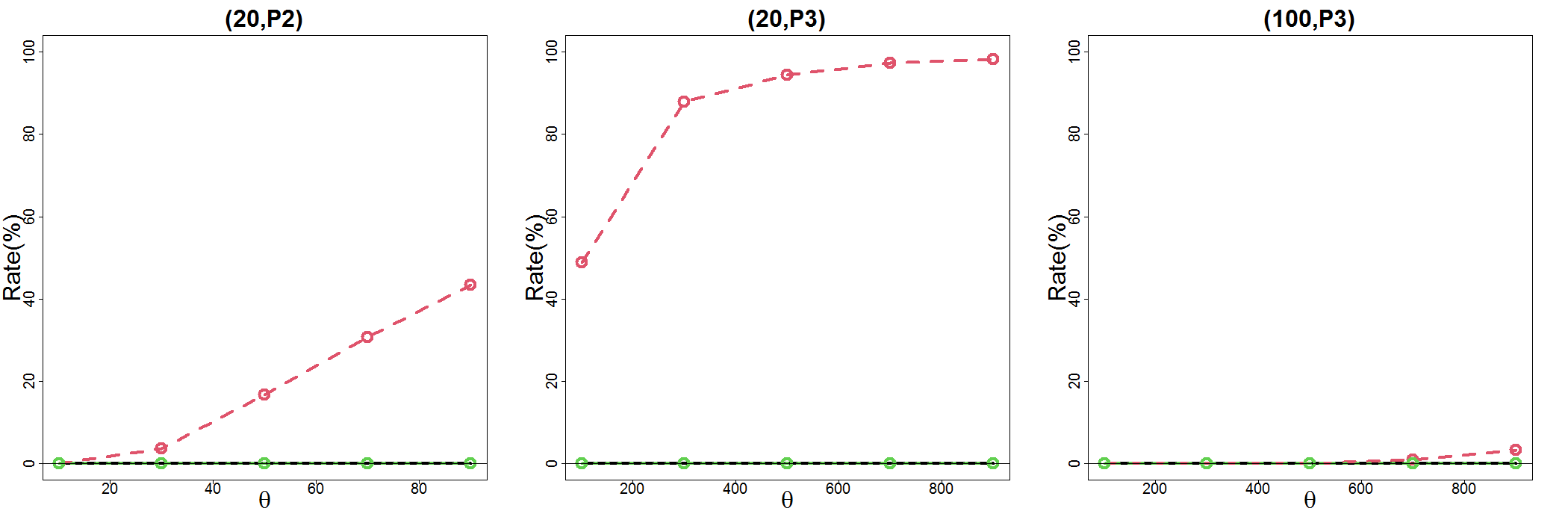}
 \caption{Rate of occurrence of negative estimates of three estimators of $R_1$ (NM: solid line; BC1: dashed line; and BC2: dotted line)} in three combinations: $(n,P)\in \{(20,P2),(20,P3), (100,P3)\}$; $x$ axis denotes values of $\theta$. 
\label{negaR.p}
\end{center}
\end{figure}

Finally, we report the rate of occurrence of negative estimates of $R_1$ in Figure \ref{negaR.p} for three cases: $(n,P)\in \{(20,P2),(20,P3), (100,P3)\}$, although we did not theoretically consider such cases in this study. We note that negative estimates did not occur in other cases. The results showed that only the estimates of $\hat R_1^{(BC1)}$ can be negative, as shown in the three figures. In particular, the results with a relatively large $\theta$ for sample size $n$ imply a high probability of $\hat R_i^{(BC1)}$ being unrealistic negative estimates.

\section{Conclusion and Discussion}

In this study, we constructed three types of asymptotic UMVUE of $R_i$: one matched the moments of UMVUE up to the second-order, whereas the others, up to the fourth order. In addition, the non-existence of a solution in parameter estimation, which is a serious practical problem, was avoided. Moreover, the common adjusted maximum likelihood estimator can be used in simultaneous inferences for $R_i$ and $R_j$ with $i\neq j$, similar to the maximum likelihood estimator.  Furthermore, we applied this methodology to assess the risk of population unique in disclosure control.

In our case, the simulation study showed that the fourth-order UMVUE outperformed the second-order one in terms of the relative bias and relative root of the mean squared error. By contrast, in cases $n<\theta$, the higher-order asymptotic UMVUE using the additive bias correction method may result in an extremely low efficiency, as indicated from the simulation study. Moreover, such an estimator may pose a high risk of being negative estimates of $R_i$ in such cases.

These simulation results imply that some theoretical differences may occur among our asymptotic UMVUEs in other asymptotic settings: $n\to \infty$ and $\theta\to \infty$ as discussed in Tsukuda (2017). In the near future, we will attempt to address these issues.

\section*{Acknowledgements}

The first author was supported in part by JSPS KAKENHI Grant 18K12758. 
The second author was supported in part by JSPS KAKENHI Grants 18H00835 and 20K03742.

\appendix
\def\thesection{Appendix \Alph{section}}

\section{Regularity conditions and related results}

We give regularity conditions for several theorems in the following.

\begin{RC}
The parameter $\theta>0$ and $i\in \mathbb Z_{>0}$ are bounded for large $n$. The sample size is  $n\geq 2$.
\end{RC}

\begin{RC}
The adjustment factor $\tilde l_{ad}(\xi)$ is in the class $\mathbb C_\xi^6$, where $\mathbb C_\xi^6$ is the set of the sixth-times differentiable functions of $\xi$ on $\mathbb R$. 
In addition, $\partial_\xi^{(j)}\tilde l_{ad}(\xi)$ does not depend on the random variable $K_n$ and are of the order $O(1)$ for $j=1,\ldots, 6$.
\end{RC}

Moreover, several calculation results are obtained, as follows:
\begin{align}
\partial_\theta R_i\equiv& \frac{\partial R_i(\xi)}{\partial \theta}=\frac{R_i}{\theta}\left[1-\sum_{j=1}^i\frac{\theta}{\theta+N-j}\right],\label{e3.0.1103}\\
\partial_{\theta}^2 R_i(\theta)=&\frac{R_i}{\theta} \left[-2\sum_{j=1}^i\frac{(N-j)}{(\theta+N-j)^2}+\sum_{j=1}^i\sum_{s\neq j}^i\frac{\theta}{(\theta+N-j)(\theta+N-s)}\right]. \label{e2.201229}
\end{align}

\begin{align}
g_\xi &\equiv g(\xi)=-\frac{\partial^2 l(\xi)}{\partial \xi^2}=\sum_{j=2}^n\frac{\theta(j-1)}{(\theta+j-1)^2},\label{e5.1103}\\
\partial_\xi g &\equiv \frac{\partial g_\xi}{\partial \xi}=\sum_{j=2}^n\frac{\theta(j-1)(j-1-\theta)}{(\theta+j-1)^3}.\label{e6.1103}\\
\partial_\xi^2 g=&\theta \partial_\theta g+\theta^2 \partial_\theta^2 g, \ 
\partial_\xi^3 g=\theta \partial_\theta g+3\theta^2 \partial_\theta^2 g+\theta^3 \partial_\theta^3 g,\label{e6.210125}\\
\partial_\xi^4 g=&\theta \partial_\theta g+7\theta^2 \partial_\theta^2 g+6\theta^3 \partial_\theta^3 g+\theta^4 \partial_\theta^4 g,\label{e7.210125} 
\end{align}
where $\xi=\log \theta$ is a natural parameter of the exponential distribution family.

\noindent In addition, under regularity condition R1 with $j=2,3,4$ and $l=1,2,3,4$, we have the following for large $n\leq N$:
\begin{align}
R_i\sim \frac{\theta}{i}=O(1), \ &\partial_\theta R_i\sim \frac{1}{i}=O(1), \ \partial_\theta^2 R_i=O(N^{-1}), \notag\\
g_\xi \sim \theta \log n=O(\log n),& \ \partial_\xi g_\xi \sim \theta \log n=O(\log n), \ \partial_{\theta}^j g_\xi=O(1), \notag\\
g_\xi-\partial_\xi^{l} g=O(1),&\ B_i(\theta)=O((\log n)^{-2}), \ \partial_\theta B_i(\theta)=O((\log n)^{-2}).\label{e1.1.18}
\end{align}
\noindent To satisfy the first and second lines above, we use \eqref{e3.0.1103}--\eqref{e6.1103}, whereas the second line and \eqref{e5.1103}--\eqref{e7.210125} are used for the last line.

Furthermore, we establish Lemma \ref{lem1} to prove the theorems related to estimator $\hat \theta_{GA}$ in \eqref{e2.210127}. Subsequently, we redefine $\mathscr {S}=\{\hat \theta: \hat \theta \in (0,C_+]\}$. 

\begin{lem}
\label{lem1}
On set $\mathscr {S}$, we have the following for large $n$ under regularity conditions R1--R2:
\begin{align*}
(i) E^{\mathscr S}[\hat \theta_{GA}-\theta]=
\frac{\theta}{g_\xi}&\left[\partial_\xi^{(1)} \tilde l_{ad}\left(1+\frac{\partial_\xi^{(2)}\tilde l_{ad}}{g_{\xi}}\right)-\frac{\partial_\xi^{(2)} \tilde l_{ad}\partial_\xi g_\xi}{2g_{\xi}^2}\right]\notag\\
&+\frac{\theta}{2g_\xi^2}\left(g_\xi-\partial_\xi g_\xi-\partial_\xi^{(2)} \tilde l_{ad}+\partial_\xi^{(3)} \tilde l_{ad}\right)+o((\log n)^{-2});\\
(ii) E^{\mathscr{S}}[(\hat \theta_{GA}-\theta)^2]=&\frac{\theta^2}{g_\xi}\left[1+\frac{(\partial_\xi \tilde l_{ad})^2+2\partial_\xi^{(2)} \tilde l_{ad}(\xi)}{g_\xi} \right]+o((\log n)^{-2}),
\end{align*}
\noindent where $E^{\mathscr{S}}[\cdot]$ is the expectation on set $\mathscr{S}$. 
\end{lem}
\noindent The proof is shown in Appendix C.3.

\section{Proof of Theorems}

\subsection{Theorems \ref{thm0}--\ref{thm2}}

For Theorem \ref{thm0}, we use a method similar to that of Das et al. (2004) and Lemma 3.3 in Hirose and Mano (2021). Under the regularity condition, the following holds on  $\mathscr S$:
\begin{align}
E^{\mathscr S}[R_{i}(\hat \theta_{ML})-R_{i}]=&E^{\mathscr S}[(\hat \theta_{ML}-\theta)]\partial_\theta R_i+\frac{1}{2}E^{\mathscr S}[(\hat \theta_{ML}-\theta)^2]\partial_\theta^2 R_i\Big|_{\theta=\theta^*},\notag\\
=&E^{\mathscr S}[(\hat \theta_{ML}-\theta)]\frac{R_i}{\theta}+o((\log n)^{-2}), \label{e7.1103}
\end{align}
where $\theta^*$ lies between $\theta$ and $\hat \theta_{ML}$. In the above, the second equality holds from Lemma \ref{lem1}, \eqref{e3.0.1103}, and \eqref{e1.1.18}.

\noindent From Lemma \ref{lem1} (i), \eqref{e3.0.1103}, \eqref{e5.1103}, and \eqref{e1.1.18}, 
the following holds for large $n$:
\begin{align*}
\eqref{e7.1103}=&\frac{R_i}{\left\{\sum_{j=2}^n\frac{(j-1)}{(\theta+j-1)^2}\right\}^2}\sum_{j=2}^n\frac{(j-1)}{(\theta+j-1)^3}+o((\log n)^{-2})=O((\log n)^{-2}).
\end{align*}

\noindent Hence, under regularity condition R1, we have the following for large $n$:

\begin{align}
|E[\hat R_i^{(NM)}-R_i]|\leq &|E^{\mathscr S}[\hat R_i^{(NM)}-R_i]|+|E^{\mathscr S_K}[\hat R_i^{(NM)}-R_i]|+|E^{\mathscr K_1}[\hat R_i^{(NM)}-R_i]|,\notag\\
\leq& \frac{R_i}{\left\{\sum_{j=2}^n\frac{(j-1)}{(\theta+j-1)^2}\right\}^2}\sum_{j=2}^n\frac{(j-1)}{(\theta+j-1)^3}+C_R\{ P(\mathscr S_K^c)+P(\mathscr K_1)\}+o((\log n)^{-2}),\notag\\
=&\frac{R_i}{\left\{\sum_{j=2}^n\frac{(j-1)}{(\theta+j-1)^2}\right\}^2}\sum_{j=2}^n\frac{(j-1)}{(\theta+j-1)^3}+o((\log n)^{-2}),\label{e2.210118}
\end{align}
\noindent where $\mathscr S_K=\mathscr S^c \cap \mathscr K_1^c$ and $C_R=\sup |\hat R_i^{NM}-R_i|$. The $C_R$ is of the order $O(1)$ for large $n$ from the definition of the estimator of $R_i$. 
For the last equality to be valid, probabilities $P(\mathscr S_K)$ and $P(\mathscr K_1)$ are of the order $o((\log n)^{-2})$ for large $n$, as a result of Lemmas \ref{lem2} and \ref{lem2.5}.

\noindent Therefore, Theorem \ref{thm0} (i) is shown from \eqref{e2.210118}, along with Theorem \ref{thm2} (i).

\vskip 1em

Next, we prove part (ii) for each theorem. From Lemma \ref{lem1} (ii), \eqref{e3.0.1103}, \eqref{e2.201229}, and \eqref{e1.1.18}, we obtain the following for large $n$:  
\begin{align*}
E^{\mathscr S}[(R_i(\hat \theta_{ML})-R_i)^2]=&E^{\mathscr S}[(\hat \theta_{ML}-\theta)^2](\partial_\theta R_i)^2+o((\log n)^{-2}),\\
=&\frac{R_i^2}{\sum_{j=2}^n\frac{\theta(j-1)}{(\theta+j-1)^2}}+o((\log n)^{-2}).
\end{align*}

\noindent Hence, Theorem \ref{thm0} (ii) is obtained in a manner similar to \eqref{e2.210118}, as follows:
\begin{align*}
E[(\hat R_i^{(NM)}-R_i)^2]= &\frac{R_i^2}{\sum_{j=2}^n\frac{\theta(j-1)}{(\theta+j-1)^2}}+o((\log n)^{-2}).
\end{align*}

For Theorem \ref{thm2} (ii), it is expressed on $\mathscr{S}$ that
\begin{align*}
E^{\mathscr S}[\{R_i(\hat \theta_{ML})-R_i-B_i(\hat \theta_{ML})\}^2]=&E^{\mathscr S}[\{R_i(\hat \theta_{ML})-R_i-(B_i(\hat \theta_{ML})-B_i)-B_i\}^2],\notag\\
=&\frac{R_i^2}{\sum_{j=2}^n\frac{\theta(j-1)}{(\theta+j-1)^2}}+o((\log n)^{-2}).
\end{align*}
\noindent We used \eqref{e1.1.18} for the second equality above. Hence, Theorem \ref{thm2} (ii) is also shown in a similar manner to \eqref{e2.210118}.

\subsection{Theorem \ref{thm3}}

For Theorem \ref{thm3}, we consider $\hat \theta_{GA}$ introduced in Section 3.2. 
The following result is obtained under regularity conditions on set $\mathscr S$:
\begin{align}
E^{\mathscr S}[R_{i}(\hat \theta_{GA})-R_{i}]=&E^{\mathscr S}[(\hat \theta_{GA}-\theta)]\partial_\theta R_i+\frac{1}{2}E^{\mathscr S}[(\hat \theta_{GA}-\theta)^2]\partial_\theta^2 R_i\Big|_{\theta=\theta^*},\notag\\
=&E^{\mathscr S}[(\hat \theta_{GA}-\theta)]\frac{R_i}{\theta}+o((\log n)^{-2}),\label{e1.210317} 
\end{align}
\noindent where $\theta^*$ lies between $\theta$ and $\hat \theta_{GA}$. In the equation above, we used Lemma \ref{lem1}, \eqref{e3.0.1103}, and \eqref{e1.1.18}.

\noindent From Lemma \ref{lem1} (i), it can be rewritten as follows for large $n$:

\begin{align*}
\eqref{e1.210317}=&\frac{R_i}{g_\xi}
\left[\partial_\xi^{(1)} \tilde l_{ad}\left(1+\frac{\partial_\xi^{(2)}\tilde l_{ad}}{g_{\xi}}\right)-\frac{\partial_\xi^{(2)} \tilde l_{ad}\partial_\xi g_\xi}{2g_{\xi}^2}\right]+\frac{R_i}{2g_\xi^2}\left(g_\xi-\partial_\xi g_\xi-\partial_\xi^{(2)} \tilde l_{ad}+\partial_\xi^{(3)} \tilde l_{ad}\right)\notag\\
&+o((\log n)^{-2}).
\end{align*}

\noindent In addition, using Lemma \ref{lem1} (ii), \eqref{e3.0.1103}, \eqref{e5.1103}, and \eqref{e1.1.18}, the following is obtained for large $n$:
\begin{align}
E^{\mathscr S}[(R_i(\hat \theta_{GA})-R_i)^2]=&
\frac{R_i^2}{g_\xi}\left[1+\frac{(\partial_\xi \tilde l_{ad})^2+2\partial_\xi^{(2)} \tilde l_{ad}(\xi)}{g_\xi} \right]+o((\log n)^{-2}),\notag
\end{align}
where $g_\xi=\sum_{j=2}^n\frac{\theta (j-1)}{(\theta+j-1)^2}$.

Hence, Theorem \ref{thm3} is calculated similarly as  \eqref{e2.210118}.

\section{Proof of Lemmas}

\subsection{Lemma \ref{lem2}}

Using \eqref{Un.dist}, probabilities $P(K_n=1)$ and $P(K_n=n)$ are calculated using the property of the unsigned Stirling number, as follows:

\begin{align*}
P(K_n=1)=&\frac{\theta}{\theta^{[n]}}s(n,1)=\frac{\theta \Gamma(\theta)(n-1)!}{\Gamma(\theta+n)}\sim \Gamma(\theta+1)n^{-\theta},\\
P(K_n=n)=&\frac{\theta^n}{\theta^{[n]}}s(n,n)=\frac{\theta}{\theta^{[n]}}
\sim \frac{\Gamma(\theta)}{\surd{2\pi}}e^n \theta^n n^{-n-\theta+1/2}.
\end{align*}
\noindent It is note that Stirling's formula was used in the calculations above.

\noindent Hence, $$P(\mathscr K)=O(n^{-\theta}\vee e^n \theta^n n^{-n-\theta+1/2})=O(n^{-\theta}).$$ This lemma is then obtained from the equalities $\lim_{n\to \infty} (\log n)^2n^{-\theta}=0$.

\subsection{Lemma \ref{lem2.5}}

Theorem 2 in Das et al. (2004) provides 
$P(\mathscr B^c)\sim O((\log n)^{-\nu/8})$ with any $\rho \in (0,1)$ and a finite positive value $\nu$.
It is note that $\mathscr B$ is a set that satisfies for large $n$, on $\mathscr B$, 
$\hat \theta \in \Theta, \partial_\theta l(\hat \theta)=0$, 
$|\hat \theta-\theta|<C_0 (\log n)^{-\rho/2}$, and 
$$\hat \theta=\theta+\frac{\partial_\theta l(\theta)}{g_\theta}+r_0,$$ where 
$|r_0|\leq C_0(\log n)^{-\rho}u_*$ with $E(u_*^\nu)$ being bounded and $C_0$ is a positive generic constant.

Let $\mathscr{S}_1$ be a set $\{\hat \theta:|\hat \theta-\theta|<2C_0\}$. Then it holds that 
$\mathscr S_1 \subseteq \mathscr S$ under regularity condition R1 from the result $\hat \theta>0$ when $1<K_n<n$. Hence, we obtain $\nu>16$,
\begin{align*}
P(\mathscr S^c)\leq P(\mathscr S_1^c)\leq P(\mathscr B^c)=o((\log n)^{-2}).  
\end{align*}

\noindent In the above, the inequality $P(\mathscr B)\leq P(\mathscr S_1)$ under regularity condition R1 is obtained from the result $|\hat \theta-\theta|<C_0 (\log n)^{-\rho/2}<2C_0$ for the case $n\geq 2$.

\subsection{Lemma \ref{lem1}}

We now establish a new lemma.

\begin{lem}
\label{lem3}
Under regularity conditions R1 and R2, for $j=3,4$, it holds that
\begin{align*}
E^{\mathscr S}[(\hat\xi_{GA}-\xi)^j]&=O((\log n)^{-2}),\\
E^{\mathscr S}[(\hat\xi_{GA}-\xi)^5]&=O((\log n)^{-3}).
\end{align*}
\end{lem}

\noindent The proof is shown from \eqref{e6.210125}--\eqref{e1.1.18} and Holder's inequality. 
In addition, the following holds for the natural parameter $\xi$ on the exponential family distribution:
\begin{align}
\mu_2=\kappa_2, \ &\mu_3=\kappa_3, \ \mu_4=\kappa_4+3\kappa_2^2, \ \mu_5=\kappa_5+10\kappa_2\kappa_3,\notag\\ 
\mu_6=&\kappa_6+15\kappa_4\kappa_2+10\kappa_3^2+15\kappa_2^3,\notag\\
\kappa_j=&\partial_\xi^{(j-2)} g_\xi, \ (j=2,\ldots,6),\label{p.lem3.1}
\end{align}
where $\mu_t$ and $\kappa_t$ are the central $t$-th moment and cumulant of $K_n$, given $n$.

In the model, we obtain the following using the natural parameter $\xi=\log \theta$:
\begin{align}
E^{\mathscr S}[\hat \theta_{GA}-\theta]=&E^{\mathscr S}\left[\sum_{i=1}^{4}\frac{(\hat \xi_{GA}-\xi)^i}{i!}\right]e^{\xi}+E^{\mathscr S}\left[\frac{(\hat \xi_{GA}-\xi)^5}{5!}\right]e^{\xi^*},\label{e1.210125}
\end{align}
where $\xi^*$ lies between $\xi$ and $\hat \xi_{GA}$ on $\mathscr S$ and 
$\hat \xi_{GA}=\log \hat \theta_{GA}$. 

Meanwhile, under the regularity conditions, it holds that
\begin{align}
-\partial_\xi l_{ad}(\xi)=&\sum_{i=1}^4\frac{(\hat \xi_{GA}-\xi)^{i}}{i!} \partial_\xi^{(i+1)} l_{ad}(\xi)+\frac{(\hat \xi_{GA}-\xi)^5}{5!}\partial_\xi^{(6)} l_{ad}(\xi)\Big|_{\xi=\xi^{**}},\notag\\
=&\sum_{i=1}^4\frac{(\hat \xi_{GA}-\xi)^{i}}{i!} \left[\partial_\xi^{(2)} l_{ad}(\xi)-\left\{\partial_\xi^{(2)} l_{ad}(\xi)-\partial_\xi^{(i+1)} l_{ad}(\xi)\right\}\right]+r_1,\label{e2.210125}
\end{align}
where $\xi^{**}$ lies between $\xi$ and $\hat \xi_{GA}$ on $\mathscr S$, and $r_1$ indicates the last terms on the right-hand side of the first line.

\noindent In Equation \eqref{e2.210125}, it follows that 
$E^{\mathscr S}[r_1]=O((\log n)^{-2})$ on $\mathscr{S}$ from Lemma \ref{lem3}, \eqref{e7.210125}, and \eqref{e1.1.18}.

\noindent Equation \eqref {e2.210125} can be rewritten as follows:
\begin{align}
\sum_{i=1}^4 \frac{(\hat \xi_{GA}-\xi)^{i}}{i!}=&\frac{1}{g_\xi}\Big[\partial_\xi^{(1)} l_{ad}+(\hat \xi_{GA}-\xi)\partial_\xi^{(2)} \tilde l_{ad}-\frac{1}{2}(\hat \xi_{GA}-\xi)^2\left(\partial_\xi^{(2)} l_{ad}-\partial_\xi^{(3)} l_{ad}\right)\notag\\
&-\frac{1}{6}(\hat \xi_{GA}-\xi)^3\left(\partial_\xi^{(2)} l_{ad}-\partial_\xi^{(4)} l_{ad}\right)
-\frac{1}{4!}(\hat \xi_{GA}-\xi)^4\left(\partial_\xi^{(2)} l_{ad}-\partial_\xi^{(5)} l_{ad}\right)\Big]+r_2,\notag\\
=&\frac{1}{g_\xi}\Big[\partial_\xi^{(1)} l_{ad}+(\hat \xi_{GA}-\xi)\partial_\xi^{(2)} \tilde l_{ad}+\frac{1}{2g_{\xi}}\left(g_\xi-\partial_\xi g_\xi-\partial_\xi^{(2)} \tilde l_{ad}+\partial_\xi^{(3)} \tilde l_{ad}\right)\Big]+r_3.\label{e1.210120}
\end{align}

\noindent Note that Lemma \ref{lem3}, \eqref{e7.210125}, and \eqref{e1.1.18} are used to show 
$E^{\mathscr S}[r_2]=o((\log n)^{-2})$, 
where $r_2$ appears on the right side of the first equality. 
On the right side of the second equality, $r_3$ is expressed as
\begin{align*}
r_3=-\frac{1}{g_\xi}\Big[&
\frac{1}{2}\left\{(\hat\xi_{GA}-\xi)^{2}-\frac{1}{g_\xi}\right\}(\partial_\xi^{(2)} l_{ad}-\partial_\xi^{(3)} l_{ad})\\
&+\frac{1}{6}(\hat\xi_{GA}-\xi)^{3}(\partial_\xi^{(2)} l_{ad}-\partial_\xi^{(4)} l_{ad})+\frac{1}{4!}(\hat\xi_{GA}-\xi)^{4}(\partial_\xi^{(2)} l_{ad}-\partial_\xi^{(5)} l_{ad})
\Big]+r_2.
\end{align*}

\noindent Moreover, the following holds for large $n$:
\begin{align}
E^{\mathscr S}[(\hat\xi_{GA}-\xi)^2]=\frac{1}{g_\xi}+o((\log n)^{-1}).\label{e1.210416}
\end{align}
\noindent Subsequently, Lemma \ref{lem3}, \eqref{e1.1.18}, and \eqref{e1.210416} yield the following result: $E^{\mathscr S}[r_3]=o((\log n)^{-2})$.

\noindent Hence, for large $n$, we obtain
\begin{align}
E^{\mathscr S}\left[\sum_{i=1}^4 \frac{(\hat \xi_{GA}-\xi)^{i}}{i!}\right]=&\frac{1}{g_\xi}\left[\partial_\xi^{(1)} \tilde l_{ad}\left(1+\frac{\partial_\xi^{(2)}\tilde l_{ad}}{g_{\xi}}\right)-\frac{\partial_\xi^{(2)} \tilde l_{ad}\partial_\xi g_\xi}{2g_{\xi}^2}\right]\notag\\
&+\frac{1}{2g_\xi^2}\left(g_\xi-\partial_\xi g_\xi-\partial_\xi^{(2)} \tilde l_{ad}+\partial_\xi^{(3)} \tilde l_{ad}\right)\notag\\
&+o((\log n)^{-2}). \label{e3.210125}
\end{align}

\noindent For the equality above to be valid, we used the following result, which is obtained from \eqref{e2.210125}:
\begin{align*}
	E^{\mathscr S}[\hat\xi_{GA}-\xi]=\frac{1}{g_\xi}\left[\partial_\xi^{(1)} \tilde l_{ad}-\frac{\partial_\xi g_\xi}{2g_\xi}\right]+o((\log n)^{-1}).
 \end{align*}

\noindent Combining Lemma \ref{lem3}, \eqref{e1.210125}, and \eqref{e3.210125}, part (i) is obtained as follows:
\begin{align*}
E^{\mathscr S}[\hat \theta_{GA}-\theta]=&
\frac{\theta}{g_\xi}\left[\partial_\xi^{(1)} \tilde l_{ad}\left(1+\frac{\partial_\xi^{(2)}\tilde l_{ad}}{g_{\xi}}\right)-\frac{\partial_\xi^{(2)} \tilde l_{ad}\partial_\xi g_\xi}{2g_{\xi}^2}\right]\notag\\
&+\frac{\theta}{2g_\xi^2}\left(g_\xi-\partial_\xi g_\xi-\partial_\xi^{(2)} \tilde l_{ad}+\partial_\xi^{(3)} \tilde l_{ad}\right)+o((\log n)^{-2}).
\end{align*}

Next, we present part (ii). On set $\mathscr{S}$, using Lemma \ref{lem3} similarly as \eqref{e1.210125}, we obtain
\begin{align}
E^{\mathscr{S}}[(\hat \theta_{GA}-\theta)^2]=&E^{\mathscr{S}}[(e^{\hat \xi_{GA}}-e^{\xi})^2],\notag\\
=&E^{\mathscr{S}}\left[\left\{\sum_{i=1}^{3}\frac{(\hat \xi_{GA}-\xi)^i}{i!}\right\}^2\right]e^{2\xi}+o((\log n)^{-2}).\label{e4.210125}
\end{align}

\noindent Furthermore, using \eqref{e1.1.18}, \eqref{e1.210120} and Lemma \ref{lem3},
\begin{align*}
\left(\sum_{i=1}^3 \frac{(\hat \xi_{GA}-\xi)^{i}}{i!}\right)^2=&\frac{1}{g_\xi^2}\left[\left(\partial_\xi^{(1)} l_{ad}\right)^2+2(\hat \xi_{GA}-\xi)\partial_\xi^{(1)} l_{ad}\partial_\xi^{(2)} \tilde l_{ad}\right]+r_4,\\
=&\frac{1}{g_\xi^2}\left[\left(\partial_\xi^{(1)} l(\xi)+\partial_\xi^{(1)} \tilde l_{ad}(\xi)\right)^2+2(\hat \xi_{GA}-\xi)\partial_\xi^{(1)} l_{ad}\partial_\xi^{(2)} \tilde l_{ad}\right]+r_4, 
\end{align*}
where $r_4$ satisfies $E[r_4]=o((\log n)^{-2})$.

\noindent Based on the result above, Equation \eqref{e4.210125} reduces to
\begin{align*}
  \eqref{e4.210125}=
  &\frac{\theta^2}{g_\xi}
  \left[1+\frac{(\partial_\xi \tilde l_{ad}(\xi))^2
      +2\partial_\xi^{(2)} \tilde l_{ad}(\xi)}{g_\xi}\right]+o((\log n)^{-2}).
\end{align*} 
Hence, Lemma \ref{lem1} is proved.

\end{document}